\documentclass[11pt]{amsart}
\usepackage[margin=2cm,bottom=1cm,includefoot]{geometry}
\usepackage{amsmath,amssymb,amsthm}
\usepackage{bbm}
\usepackage{enumerate}
\usepackage{hyperref}
\usepackage[all]{hypcap}

\newtheorem{thm}{Theorem}
\newtheorem{prop}[thm]{Proposition}
\newtheorem{cor}[thm]{Corollary}
\theoremstyle{definition}
\newtheorem{rem}[thm]{Remark}

\numberwithin{equation}{section}

\newcommand{\N}{\mathbb{N}}
\newcommand{\E}{\mathbb{E}}
\newcommand{\p}{\mathbb{P}}
\newcommand{\1}{\mathbbm{1}}
\DeclareMathOperator{\e}{e}

\newcommand{\op}{\operatorname}
\newcommand{\lne}{<}
\newcommand{\gne}{>}
\newcommand{\abs}{|}

\begin{document}
\title[Inversions in Random Permutations Under the Ewens Sampling Distribution]{Inversions in Random Permutations \linebreak Under the Ewens Sampling Distribution \linebreak With and Without a Prescribed Number of Fixed Points}
\author{Ross G. Pinsky, Dominic T. Schickentanz \medskip\\ November 16, 2025 }
\thanks{The second named author has been supported in part by a Technion fellowship.}


\maketitle

\begin{abstract}
	In the first part of the paper, we study the inversion statistic of random permutations under the family $(\mathbb{P}_\theta^{(n)})_{\theta \ge 0}$ of Ewens sampling distributions on $S_n$. We obtain a rather simple exact formula for the expected number of inversions under $\mathbb{P}_\theta^{(n)}$. In particular, we show that this expected number of inversions is decreasing in the tilting parameter $\theta$ for any $n$ and that it is convex in $\theta$ for $n \not \in \{3,4\}$ only. Furthermore, we derive an exact formula for the probability that a specific pair of indices $(i,j) \in \{1,\dots,n\}^2$ is inverted and show that this probability is decreasing in $\theta$ if and only if $|j-i| \ge 2$ holds. We also exhibit the asymptotic behavior of these quantities as $n \to \infty$ and $\theta \to \infty$.
	
	In the second part of our paper, we analyze the inversion statistic of random permutations under~$(\mathbb{P}_\theta^{(n)})_{\theta > 0}$ conditioned on having a prescribed number of fixed points. Again, we obtain exact formulas for the expected number of inversions and for the probability that a specific pair of indices is inverted. Since, as expected, the resulting formulas are rather complicated, we focus on the asymptotic behavior of these quantities as $n \to \infty$, $\theta \to \infty$ and $\theta \to 0$.
\end{abstract}

\section{Introduction and Statement of Results}
In this paper, we study the inversion statistic of random permutations under the family of Ewens sampling distributions, the distributions obtained by exponential tilting via the total-number-of-cycles statistic. We consider both the unconditioned Ewens sampling distributions and the Ewens sampling distributions conditioned on the permutation having a prescribed number of fixed points. Surprisingly, it does not seem that the inversion statistic has been studied even under the unconditioned Ewens sampling distributions.
Our results in this regard may be thought of as the obverse of the results in~\cite{GP18} concerning cycle statistics under the Mallows distributions, the distributions obtained by exponential tilting via the inversion statistic.

Our results dealing with permutations under Ewens sampling distributions conditioned on having a prescribed number of fixed points generalize recent results of~\cite{Pin25}. In that paper, one of the present authors studied the inversion statistic under uniformly distributed random permutations conditioned on having a prescribed number of fixed points. The proofs in that paper exploited the so-called Chinese restaurant construction of a uniformly distributed random permutation. As is well-known, that construction can be modified to produce a random permutation with a Ewens sampling distribution~\cite{Ald85,Pit06}. We use this construction in our proofs.

Throughout this paper, let~$n \in \N$ with~$n \ge 3$ (unless explicitly indicated otherwise) and~$\theta \ge 0$. As usual, we denote the set of all permutations of~$[n]:=\{1,\dots,n\}$ by~$S_n$. Let $$N(\pi):= \max\{\abs A\abs: A \subseteq [n],\, \pi^m(i) \neq j \text{ for all } i,j \in A \text{ with } i \neq j \text{ and all } m \in \N\}, \quad \pi \in S_n,$$ be the number of cycles in~$\pi$. Then the Ewens sampling distribution on~$S_n$ with parameter~$\theta$ is the probability measure~$\p_\theta^{(n)}$ on~$(S_n,\mathcal{P}(S_n))$ defined by~$$\p_\theta^{(n)}(\{\pi\}) := \frac{\theta^{N(\pi)}}{\theta^{(n)}}, \quad \pi \in S_n,$$ where the normalization constant~$\theta^{(n)}$ is the rising factorial defined by
$\theta^{(n)}= \prod_{k=0}^{n-1} (\theta+k)$. (In the special case~$\theta=0$, we understand~$\frac{0^{N(\pi)}}{0^{(n)}}$ to be~$\frac{\1_{N(\pi)=1}}{(n-1)!}$.)

As~$\theta$ increases, permutations with more (and hence shorter) cycles become more likely under~$\p_\theta^{(n)}$, interpolating between the uniform distribution on the set of cyclic permutations for~$\theta=0$ and the point measure on the identity permutation~$\op{id}_{[n]}$ as~$\theta \to \infty$. Of course,~$\p_1^{(n)}$ is the uniform distribution on~$S_n$. For later use, we recall one well-known property \cite{ABT03,Pin14}: As~$n\to\infty$, the number of cycles of length~$m \in \N$ of a $\mathbb{P}_\theta^{(n)}$-distributed random permutation converges in distribution to the Poisson distribution with parameter~$\frac\theta m$. It is convenient to let~$\Sigma_n: S_n \to S_n$ denote the identity in order to have a generic~$\p_\theta^{(n)}$-distributed random variable under~$\p_\theta^{(n)}$ at hand.

Originating from population genetics~\cite{Ewe72}, Ewens sampling has been applied in several non-mathe\-matical fields, ranging from ecology to physics. For instance, it plays a crucial role in the unified neutral theory of biodiversity. At the same time, it appears in various areas of pure mathematics, ranging from algebra and number theory to probability. We refer to~\cite{JKB97} and~\cite{Cra16} for more insights on the ubiquity of Ewens sampling.

Throughout this paper, let~$i,j \in [n]^2$ with~$i \lne j$.
We recall that the pair~$(i,j)$ is called an inversion under~$\pi \in S_n$ if~$\pi^{-1}(i) \gne \pi^{-1}(j)$ holds, that is, if~$i$ appears to the right of~$j$ in the one-line notation of~$\pi$. We denote by $$\op{Inv}(\pi):= \{(i,j) \in [n]^2: i \lne j,\ (i,j) \text{ is an inversion}\}, \quad \pi \in S_n,$$ the set of all inversions in the permutation~$\pi$. In the literature, there also is an alternative convention defining inversions through images instead of preimages. Noting~$\Sigma_n^{-1} \sim \Sigma_n$ under~$\p_\theta^{(n)}$, the difference between the two conventions is irrelevant in our probabilistic framework.

\subsection{Inversions Under (Unconditioned) Ewens Sampling}
\label{sec1}
By symmetry, under the uniform distribution~$\p_1^{(n)}$, one has~$\E_1^{(n)} \abs\op{Inv}(\Sigma_n)\abs =\frac14 n(n-1)$. Would one expect~$\E_\theta^{(n)} \abs\op{Inv}(\Sigma_n)\abs$ to be larger or smaller than~$\E_1^{(n)} \abs\op{Inv}(\Sigma_n)\abs$? Of course, this should depend on the value of~$\theta$.
On the one hand, recall from above that the number of fixed points in a $\p_\theta^{(n)}$-distributed permutation converges in distribution to the Poisson distribution with parameter~$\theta$, as~$n \to \infty$. The Poisson distributions are stochastically increasing with respect to their parameter~$\theta$. Having a lot of fixed points decreases the number of inversions since each pair of fixed points automatically does not form an inversion. For large~$n$, this suggests that~$\E_\theta^{(n)} \abs\op{Inv}(\Sigma_n)\abs$ should be larger than~$\E_1^{(n)} \abs\op{Inv}(\Sigma_n)\abs$ for~$\theta\in(0,1)$ and vice versa for~$\theta>1$. On the other hand, recall from above that the number of two-cycles in a~$\p_\theta^{(n)}$-distributed permutation converges in distribution to the Poisson distribution with parameter~$\frac\theta2$, as~$n \to \infty$.	Having a lot of two-cycles increases the number of inversions since each pair in a two-cycle automatically forms an inversion. This would suggest the opposite conclusion from the one suggested above.

\pagebreak
Here is our main result for unconditioned Ewens sampling distributions.
\begin{thm}
	\label{thm:uncond}
	\begin{enumerate}[$(a)$]
		\item The probability that~$(i,j)$ is an inversion is given by
		\begin{align}
			\label{eqn:P}
			\p_\theta^{(n)}((i,j) \in \op{Inv}(\Sigma_n))
			={}& \frac{n(n-2(j-i)+1)}{2 (\theta+n-1)} - \frac{(n-1)(n-2(j-i))}{2 (\theta+n-2)}.
		\end{align}
		As a function of~$\theta \ge 0$, this probability is (strictly) decreasing if and only if~$j-i \ge 2$ holds.
		\item The expected number of inversions under~$\p_\theta^{(n)}$ is given by
		\begin{align}
			\label{eqn:E}
			\E_\theta^{(n)}\abs\op{Inv}(\Sigma_n)\abs
			={}& \frac{(n+1)n^2(n-1)}{12 (\theta+n-1)} - \frac{n(n-1)^2(n-2)}{12 (\theta+n-2)},
		\end{align}
		which is strictly decreasing in~$\theta \ge 0$. It is a (strictly) convex function of~$\theta$ if and only if~$n \not\in \{3,4\}$ holds.
	\end{enumerate}
\end{thm}

The theorem immediately yields the following asymptotic behavior as~$n \to \infty$ for the quantities on the left-hand sides of~\eqref{eqn:P} and~\eqref{eqn:E}.
\begin{cor}
	\label{cor:asympExpan}
	\begin{enumerate}[$(a)$]
		\item Let~$(i_n)_n,(j_n)_n \subseteq \N$ be sequences with~$i_n \lne j_n \le n$ for all~$n \in \N$. For any fixed~$\theta \ge 0$, we have
		\begin{align*}
			&\p_\theta^{(n)}((i_n,j_n) \in \op{Inv}(\Sigma_n))\\
			&\qquad= \frac{1}{2} - \frac{(\theta-1)(j_n-i_n)}{n^2} - \frac{(\theta-1)(\theta-2)}{2n^2} + \frac{(2\theta-3)(\theta-1)(j_n-i_n)}{n^3} + O\left(\frac{1}{n^3}\right), \quad n \to \infty.
		\end{align*}
		\item For any fixed~$\theta \ge 0$, we have
		\begin{align}
			\label{eqn:Easymp}
			\E_\theta^{(n)}\abs\op{Inv}(\Sigma_n)\abs
			={}& \frac{n(n-1)}{4} - \frac{\theta-1}{6} n + \frac{\theta(\theta-1)}{12} + O\left(\frac{1}{n^2}\right), \quad n \to \infty.
		\end{align}
	\end{enumerate}
\end{cor}
\noindent Note that there is no term of order~$\frac{1}{n}$ in the expansion of the expected number of inversions.

\begin{rem}
	\label{rem:1}
	In~\cite{Pin25}, it was shown that the expected number of inversions in a random permutation in~$S_n$ distributed according to the uniform distribution conditioned on the permutation having~$m\in\mathbb{N}_0$ fixed points is given by $$\frac{n(n-1)}4-\frac{m-1}6 n-\frac{m^2-m-1}{12}+O\left(\frac1{(n-m-1)!}\right), \quad n \to \infty.$$ From the Chinese restaurant construction, it follows easily that the expected number of fixed points of~$\Sigma_n$ under~$\p_\theta^{(n)}$ is equal to~$\theta$. Note that if one sets~$\theta=m$ in~\eqref{eqn:Easymp}, one finds that the asymptotic expansions in powers of~$n$ for the expected number of inversions under the two measures coincide up to the~$O(1)$-term, where they differ.
\end{rem}

We may also allow~$\theta$ to depend on~$n$. Given~$c \gne 0$, a simple analysis of~\eqref{eqn:E} yields
\begin{align*}
	\E_{c n^\alpha}^{(n)} \abs\op{Inv}(\Sigma_n)\abs \sim
	\begin{cases}
		\dfrac{1}{3c} n^{3-\alpha}, \quad& \alpha \gne 1,\medskip\\
		\dfrac{4c+3}{12(c+1)^2} n^2, \quad& \alpha=1,\medskip\\
		\dfrac{1}{4} n^2, \quad& \alpha \lne 1,\\
	\end{cases}
	\qquad\quad n \to \infty.
\end{align*}
In particular, we have~$\lim_{n \to \infty} \E_{c n^\alpha}^{(n)} \abs\op{Inv}(\Sigma_n)\abs = 0$ if and only if~$\alpha \gne 3$. However, one can in fact show that~$\lim_{n \to \infty} \p_{c n^\alpha}^{(n)}(\Sigma_n=\op{id}_{[n]}) = 1$ is equivalent to~$\alpha \gne 2$. This follows from the the definition of the Ewens sampling distribution by taking the logarithm and applying the mean value theorem.

For fixed~$n$, Theorem~\ref{thm:uncond} extends formulas, previously known in the elementary special cases~$\theta\in\{0,1\}$ only, to the whole parameter range~$\theta \ge 0$.
For~$\theta=1$, we recover~$$\p_1^{(n)}((i,j) \in \op{Inv}(\Sigma_n))=\frac{1}{2} \qquad \text{and} \qquad \E_1^{(n)}\abs\op{Inv}(\Sigma_n)\abs = \frac{n(n-1)}{4}.$$ For~$\theta=0$, we recover $$\p_0^{(n)}((i,j) \in \op{Inv}(\Sigma_n)) = \frac{1}{2} + \frac{j-i-1}{(n-1)(n-2)} \qquad \text{and} \qquad \E_0^{(n)}\abs\op{Inv}(\Sigma_n)\abs = \frac{n(3n-1)}{12},$$ the latter being equivalent to OEIS sequence A227404, which counts the total number of inversions in all single cycle permutations of order~$n$.

At the other end of the parameter range, Theorem~\ref{thm:uncond} implies the previously unknown asymptotics
\begin{equation}
	\label{eqn:asmpP}
	\lim_{\theta \to \infty} \theta \p_\theta^{(n)} ((i,j) \in \op{Inv}(\Sigma_n)) = n-(j-i)
\end{equation}
and
\begin{equation}
	\label{eqn:asmpE}
	\lim_{\theta \to \infty} \theta \E_\theta^{(n)}\abs\op{Inv}(\Sigma_n)\abs = \frac{1}{4} \binom{2n}{3}.
\end{equation}

The proofs of the results presented in this subsection can be found in Section~\ref{secP1}. At the end of that section, we also provide a rather short and elementary proof of~\eqref{eqn:asmpP} and~\eqref{eqn:asmpE}, without relying on Theorem~\ref{thm:uncond}.


\subsection{Inversions Under Ewens Sampling With Prescribed Number of Fixed Points}
\label{sec2}
We define~$$\op{Fix}(\pi) := \{k \in [n]: \pi(k)=k\}, \quad \pi \in S_n,$$ to be the set of fixed points of~$\pi$ and set $$D_{n,m}:=\{\abs\op{Fix}(\Sigma_n)\abs = m\} = \{\Sigma_n \text{ has exactly } m \text{ fixed points}\}, \quad m \in \mathbb{Z},\, n \in \N.$$
In addition to our general assumption~$n \ge 3$, we henceforth assume~$\theta \gne 0$ and~$m \in \{0,\dots,n-2\}$ (unless explicitly indicated otherwise) to exclude trivial cases and, in particular, to ensure~$\p_\theta^{(n)}(D_{n,m}) \in (0,1)$. Our aim is to analyze the inversion statistic under the conditional law~$\p_\theta^{(n)}(\,\cdot\,\abs D_{n,m})$. An important special case is given by~$m=0$, with~$D_{n,0}$ being the event that~$\Sigma_n$ is a derangement.

As a first step, we derive a formula for~$\p_\theta^{(n)}(D_{n,m})$. While the result seems to be rather well-known for derangements, we did not find a reference for the general case.	For convenience, let us set~$\p_\theta^{(0)}(D_{0,0}):=1$ (without assigning any meaning to~$D_{0,0}$ or~$\p_\theta^{(0)}$) and~$\theta^{(0)}:=1$.
\begin{prop}
	\label{prop:pFix}
	We have
	\begin{align*}
	\p_\theta^{(n)}(D_{n,m})
	={}& \binom{n}{m} \frac{\theta^m\theta^{(n-m)}}{\theta^{(n)}}\p_\theta^{(n-m)}(D_{n-m,0})\\
	={}& \frac{n!\theta^m}{m!\theta^{(n)}} \sum_{k=0}^{n-m} \frac{(-\theta)^k \theta^{(n-m-k)}}{k!(n-m-k)!}, \qquad m \in \{0,\dots,n\},\, n \in \N.
	\end{align*}
\end{prop}
\noindent This proposition is complemented by the trivial observation
\begin{equation}
	\label{eqn:pFixTriv}
	\p_\theta^{(n)}(D_{n,m})=0, \quad m \in \mathbb{Z} \setminus \{0,\dots,n-2,n\},\, n \in \N.
\end{equation}

Under~$\p_\theta^{(n)}(\,\cdot\,\abs D_{n,m})$, the probability that the pair~$(i,j)$ is inverted and the expected number of inversions, unsurprisingly, are significantly more complicated than their unconditional analogs derived in Theorem~\ref{thm:uncond}. For this reason, we postpone the exact formulas for these quantities to Proposition~\ref{prop:exact} in the proof section. Instead, our main results focus on the behavior of these quantities as~$n \to \infty$,~$\theta \to \infty$ or~$\theta \to 0$.

We start with the asymptotics as~$n \to \infty$.
\begin{thm}
	\label{thm:condN}
	\begin{enumerate}[$(a)$]
		\item Let~$(i_n)_n,(j_n)_n \subseteq \N$ be sequences with~$i_n \lne j_n \le n$ for all~$n \in \N$. For any fixed~$\theta \gne 0$ and~$m \in \N_0$, we have
		\begin{align*}
			\p_\theta^{(n)}\big((i_n,j_j) \in \op{Inv}(\Sigma_n)\big\abs D_{n,m}\big)
			= \frac{1}{2} -& \frac{(m-1)(j_n-i_n)}{n^2} - \frac{m^2-3m-\theta+2}{2n^2}\\ &+ \frac{(m^2-4m-\theta+3)(j_n-i_n)}{n^3}
			+O\left(\frac{1}{n^3}\right), \qquad n \to \infty.
		\end{align*}
		\item For any fixed~$\theta \gne 0$ and~$m \in \N_0$, the expected number of inversions under~$\p_\theta^{(n)}(\,\cdot\,\abs D_{n,m})$ satisfies
		\begin{align}
			\label{eqn:condE}
			\E_\theta^{(n)}\big[\abs\op{Inv}(\Sigma_n)\abs \big\abs D_{n,m}\big]
			={}& \frac{n(n-1)}{4} - \frac{m-1}{6} n - \frac{m^2-m-\theta}{12} + O\left(\frac{1}{n}\right), \quad n \to \infty.
		\end{align}
	\end{enumerate}
\end{thm}

\noindent We emphasize that the $O(n)$-term in~\eqref{eqn:condE} does not depend on~$\theta$ and hence is universal for the whole family of Ewens sampling distributions.

\begin{rem}
Similar to Remark~\ref{rem:1}, a comparison with Corollary~\ref{cor:asympExpan} seems natural. As mentioned in the introduction,
we have~$\lim_{n \to \infty} \p_\theta^{(n)}(D_{n,m})=\p(Y_\theta=m)$ for a Poisson random variable~$Y_\theta$ with parameter~$\theta$. Formally setting~$m=Y_\theta$ in~\eqref{eqn:condE} and then taking expectations, we observe that the resulting expansion
\begin{align*}
	\frac{n(n-1)}{4} - \frac{\E Y_\theta -1}{6} n - \frac{\E Y_\theta^2 - \E Y_\theta -\theta}{12} + O\left(\frac{1}{n}\right)
	= \frac{n(n-1)}{4} - \frac{\theta-1}{6} n - \frac{\theta(\theta-1)}{12} + O\left(\frac{1}{n}\right)
\end{align*}
coincides with the expansion of
\begin{align*}
	\E_\theta^{(n)} \abs\op{Inv}(\Sigma_n)\abs = \sum_{m \in \N_0} \E_\theta^{(n)}\big[\abs\op{Inv}(\Sigma_n)\abs \big\abs D_{n,m}\big] \p_\theta^{(n)}(D_{n,m})
\end{align*}
provided in~\eqref{eqn:Easymp} up to the $O(1)$-term, where they differ in their sign only.
\end{rem}

\begin{rem}
For~$\theta=1$, corresponding to uniform sampling,~\eqref{eqn:condE} is consistent with the asymptotic expansion derived in~\cite{Pin25}. In that paper, the remainder term~$O\left(\frac{1}{n}\right)$ was shown to actually be~$O\big(\frac{1}{(n-m-1)!}\big)$ for~$\theta=1$. We strongly believe that such an improvement is not possible for~$\theta \neq 1$.
\end{rem}

Let us now fix~$n$ and consider~$\theta \to \infty$. Then the probability that the pair~$(i,j)$ is inverted and the expected number of inversions converge to strictly positive limits. Notably, these limits depend crucially on the parity of~$n-m$, the number of non-fixed points.
\begin{thm}
\label{thm:condThetaInf}
	\begin{enumerate}[$(a)$]
		\item We have
		\begin{align*}
			\lim_{\theta \to \infty} \p_\theta^{(n)}&\big((i,j) \in \op{Inv}(\Sigma_n)\big\abs D_{n,m}\big)\\
			&{}=
			\begin{cases}
				\frac{n-m}{n(n-1)}\left(\frac{n-m}{2}+\frac{m(n-1-(j-i))}{n-2}\right)	, \quad& \text{if } n-m \text{ is even},\medskip\\
				\frac{1}{n(n-1)}\left(\frac{(n-m+3)(n-m-3)}{2} + \frac{m(n-m)(n-1-(j-i))+3(n+(j-i)-3)}{n-2}\right), \quad& \text{if } n-m \text{ is odd}.
			\end{cases}
		\end{align*}
		\item The expected number of inversions under~$\p_\theta^{(n)}(\,\cdot\,\abs D_{n,m})$ satisfies
		\begin{align*}
			\lim_{\theta \to \infty} \E_\theta^{(n)}\big[\abs\op{Inv}(\Sigma_n)\abs \big\abs D_{n,m}\big]
			&{}=
			\begin{cases}
				\dfrac{(n-m)(3n+m)}{12}, \quad& \text{if } n-m \text{ is even}, \medskip\\
				\dfrac{(n-m+3)(n-m-3)}{4} + \dfrac{m(n-m)}{3} +2, \quad& \text{if } n-m \text{ is odd}.
			\end{cases}
		\end{align*}
	\end{enumerate}
\end{thm}

Noting~$\lim_{\theta \to \infty} \p_\theta^{(n)}(D_{n,n})=1$, the positivity of the above limits does not contradict~\eqref{eqn:asmpP} and~\eqref{eqn:asmpE}, which state the asymptotics as~$\theta \to \infty$ for the unconditional problem. Rather, Theorem~\ref{thm:condThetaInf} and Proposition~\ref{prop:pFix} imply
\begin{align*}
	\E_\theta^{(n)} \abs\op{Inv}(\Sigma_n)\abs
	={}& \sum_{m=0}^{n-2} \E_\theta^{(n)}\big[\abs\op{Inv}(\Sigma_n)\abs \big\abs D_{n,m}\big] \p_\theta^{(n)}(D_{n,m})\\
	\sim{}& \E_\theta^{(n)}\big[\abs\op{Inv}(\Sigma_n)\abs \big\abs D_{n,n-2}\big] \p_\theta^{(n)}(D_{n,n-2})
	\sim \frac{2n-1}{3}\cdot \frac{n(n-1)}{2\theta} = \frac{1}{4\theta} \binom{2n}{3}, \qquad \theta \to \infty,
\end{align*}
recovering~\eqref{eqn:asmpE}. Likewise,~\eqref{eqn:asmpP} may be recovered as well. In contrast to~\eqref{eqn:asmpP}, various different scenarios in terms of the (modified) Chinese restaurant construction contribute to the limit in Theorem~\ref{thm:condThetaInf}$(a)$, particularly if~$n-m$ is odd. Together with the resulting complex dependence on~$(n,m)$, this makes the existence of a short elementary proof of the above theorem, resembling the one we provide for~\eqref{eqn:asmpP} and~\eqref{eqn:asmpE}, seem very unlikely.

The dependence on the parity in Theorem~\ref{thm:condThetaInf} can be explained as follows: As~$\theta \to \infty$, the probability measure~$\p_\theta^{(n)}(\,\cdot\,\abs D_{n,m})$ concentrates on those permutations with exactly~$m$ fixed points which have as short and as many cycles as possible. If~$n-m$ is even, the ideal setting is given by~$\frac{n-m}{2}$ two-cycles in additional to the~$m$ fixed points. If~$n-m$ is odd, the ideal setting is given by exactly~$1$ three-cycle and~$\frac{n-m-3}{2}$ two-cycles in additional to the~$m$ fixed points.

Our final result deals with the limiting behavior as~$\theta \to 0$. Here the parity of~$n-m$ is not relevant. Instead, the cases~$m=n-3$ and~$m=n-2$ are distinguished.
\begin{prop}
	\label{prop:condTheta0}
	\begin{enumerate}[$(a)$]
		\item We have
		\begin{align*}
			\lim_{\theta \to 0} \p_\theta^{(n)} \big((i,j) \in \op{Inv}(\Sigma_n)\big\abs D_{n,m}\big)
			{}=
			\begin{cases}
			\tfrac{n-m}{n(n-1)} \left(\tfrac{n-m-3}{2} + \tfrac{m(n-1-(j-i))
			+ (n+j-i-3)\1_{m \neq n-3}}{n-2} \right), \quad& m \lne n-2,\medskip\\
			\tfrac{2(n-(j-i))}{n(n-1)}, \quad& m=n-2.
			\end{cases}
		\end{align*}
		\item The expected number of inversions under~$\p_\theta^{(n)}(\,\cdot\,\abs D_{n,m})$ satisfies
		\begin{align*}
			\lim_{\theta \to 0} \E_\theta^{(n)}\big[\abs\op{Inv}(\Sigma_n)\abs \big\abs D_{n,m}\big]
			={}&
			\begin{cases}
			(n-m)\left(\dfrac{n-m-3}{4} + \dfrac{m+2\1_{m \neq n-3}}{3}\right), \quad& m \lne n-2,\medskip\\
			\dfrac{2n-1}{3}, \quad& m = n-2.
			\end{cases}
		\end{align*}
	\end{enumerate}
\end{prop}

In the special case~$m=0$, corresponding to derangements, these limits are consistent with the formulas provided in Subsection~\ref{sec1}. More precisely, observing~$\p_0^{(n)}(D_{n,0})=1$, we get
$$\lim_{\theta \to 0} \E_\theta^{(n)}\big[\abs\op{Inv}(\Sigma_n)\abs \big\abs D_{n,0}\big] = \frac{n(3n-1)}{12} = \E_0^{(n)} \abs \op{Inv}(\Sigma_n)\abs = \E_0^{(n)}\big[\abs\op{Inv}(\Sigma_n)\abs \big\abs D_{n,0}\big],$$ and likewise for the probability that~$(i,j)$ is inverted.

The proofs of the results presented in this subsection can be found in Section~\ref{secP2}.

\section{Proofs of the Results in Subsection~\ref{sec1}}
\label{secP1}
The proof of Theorem~\ref{thm:uncond} relies upon a version of the Chinese restaurant construction of random permutations. A similar approach was developed recently in~\cite{Pin25} to analyze inversion statistics of uniformly random permutations with a prescribed number of fixed points.

For the reader's convenience, let us briefly recall the standard Chinese restaurant construction of random permutations generated by Ewens sampling. Consider a restaurant with an unlimited number of round tables, each of which is able to accommodate an unlimited number of people. Arriving persons take their seats according to the following iterative scheme of Markovian type: The first person simply sits at a table. Now suppose that~$n \in \N$ persons have already been seated. Then person~$n+1$ chooses to sit to the left of any particular already seated person with probability~$\frac{1}{\theta+n}$, and chooses to sit at an empty table with probability~$\frac{\theta}{\theta+n}$.
Once exactly~$n$ persons, numbered from~$1$ to~$n$, are seated, we define a random permutation~$\Sigma_n \in S_n$ by~$\Sigma_n(i):=j$, for~$i,j \in [n]$, if person~$j$ is seated to the left of person~$i$ at one of the tables. (Persons sitting alone at a table are considered their own seatmates.) One can easily verify~$\Sigma_n \sim \p_\theta^{(n)}$ under the probability measure corresponding to the $n$-step procedure described above \cite{Ald85,Pit06}. A notable feature of this construction is its consistency: Deleting~$n$ from the cycle representation of~$\Sigma_n$ yields the cycle representation of~$\Sigma_{n-1}$.

\begin{proof}[Proof of Theorem~\ref{thm:uncond}]
	$(a)$ Let~$\Sigma_n$ and~$\Sigma_{n-1}$ be the (consistent) random permutations of~$[n]$ and~$[n] \setminus \{j\}$, respectively, arising from a Chinese restaurant construction as described above, modified such that the persons	(numbered from~$1$ to~$n$) arrive in the order $$1,\dots,i-1,i+1,\dots,j-1,j+1,\dots,n,i,j.$$ We observe~$\Sigma_n \sim \p_\theta^{(n)}$ under the probability measure encoding the random seating choices of the~$n$ persons. Slightly abusing notation, we denote this probability measure by~$\p_\theta^{(n)}$ as well.
	In what follows, we split the event~$$\{(i,j) \in \op{Inv}(\Sigma_n)\}=\{\Sigma_n^{-1}(i) \gne \Sigma_n^{-1}(j)\}$$ into several sub-events and analyze them using the modified Chinese restaurant construction.
	
	First assume that person~$i$ chooses to sit to the left of~$k \in [n] \setminus \{i,j\}$ and that person~$j$ chooses to sit neither at a new table nor to the right of~$i$.
	Then~$i$ remains the seatmate to the left of~$k$ after~$j$ is seated (i.e.,~$\Sigma_n^{-1}(i)=k$). Hence~$(i,j)$ is an inversion in this scenario if and only if person~$j$ sits down to the left of a person with a label smaller than~$k$:
	\begin{align*}
		&\p_\theta^{(n)}((i,j) \in \op{Inv}(\Sigma_n), \Sigma_{n-1}(k)=i, \Sigma_n(j) \not \in \{i,j\})\\
		&\qquad= \p_\theta^{(n)}(\Sigma_n^{-1}(j)\lne k, \Sigma_n^{-1}(j) \neq j, \Sigma_{n-1}(k)=i)\\
		&\qquad= \frac{k-1-\1_{k\gne j}}{\theta+n-1} \cdot \frac{1}{\theta+n-2}.
	\end{align*}
	Substituting~$l=k-1$ for~$k\lne i$ and~$l=k-2$ for~$k\gne i$, separating the additional~$+1$ from each of the~$j-i-1$ terms with indices~$l \in \{i-1,\dots,j-3\}$ and inserting the missing summand~$l=j-2$, we deduce
	\begin{align*}
		&\p_\theta^{(n)}((i,j) \in \op{Inv}(\Sigma_n),\Sigma_{n-1}(i)\neq i, \Sigma_n(j) \not\in \{i,j\})\\
		&\qquad= \sum_{k \in [n] \setminus \{i,j\}} \p_\theta^{(n)}((i,j) \in \op{Inv}(\Sigma_n), \Sigma_{n-1}(k)=i, \Sigma_n(j) \not \in \{i,j\})\\
		&\qquad= \sum_{\substack{k=1\\k\neq i}}^{j-1} \frac{k-1}{(\theta+n-1)(\theta+n-2)} + \sum_{k=j+1}^n \frac{k-2}{(\theta+n-1)(\theta+n-2)}\\
		&\qquad= \frac{\sum_{l=0}^{n-2} l +(j-i-1) -(j-2)}{(\theta+n-1)(\theta+n-2)}\\
		&\qquad= \frac{\binom{n-1}{2} +(j-i-1) -(j-2)}{(\theta+n-1)(\theta+n-2)}.
	\end{align*}
	Next consider the case that person~$i$ chooses not to sit at a new table while person~$j$ chooses to sit at a new table. Then the right neighbor of~$i$ remains the same after~$j$ is seated (i.e.,~$\Sigma_{n-1}^{-1}(i) = \Sigma_n^{-1}(i)$) and, recalling our modified version of the Chinese restaurant construction, we get
	\begin{align*}
		&\p_\theta^{(n)}((i,j) \in \op{Inv}(\Sigma_n), \Sigma_{n-1}(i)\neq i, \Sigma_n(j)=j)\\
		&\qquad= \p_\theta^{(n)}(\Sigma_{n-1}^{-1}(i) \gne j, \Sigma_n(j)=j)\\
		&\qquad= \frac{n-j}{\theta+n-2} \cdot \frac{\theta}{\theta+n-1}.
	\end{align*}
	Now assume that person~$i$ chooses not to sit at a new table and that person~$j$ chooses to sit to the right of~$i$. Then the original right neighbor of~$i$ becomes the right neighbor of~$j$ (i.e.,~$\Sigma_{n-1}^{-1}(i) = \Sigma_n^{-1}(j)$) and we obtain
	\begin{align*}
		&\p_\theta^{(n)}((i,j) \in \op{Inv}(\Sigma_n), \Sigma_{n-1}(i)\neq i, \Sigma_n(j)=i)\\
		&\qquad= \p_\theta^{(n)}(j \gne \Sigma_{n-1}^{-1}(i), \Sigma_{n-1}^{-1}(i)\neq i, \Sigma_n(j)=i)\\
		&\qquad= \frac{j-2}{\theta+n-2} \cdot \frac{1}{\theta+n-1}.
	\end{align*}
	Finally, consider the case that person~$i$ chooses to sit at a new table. Then $(i,j)$ becomes an inversion if and only if person~$j$ sits down at the same table or to the left of a person with a label smaller than~$i$:
	\begin{align*}
		\p_\theta^{(n)}((i,j) \in \op{Inv}(\Sigma_n), \Sigma_{n-1}(i)=i)
		={}& \p_\theta^{(n)}(\Sigma_n^{-1}(j) \le i, \Sigma_{n-1}(i)=i)
		= \frac{i}{\theta+n-1}\cdot\frac{\theta}{\theta+n-2}.
	\end{align*}
	Adding the last four equations, we obtain
	\begin{align}
		\label{eqn:Pfirst}
		\p_\theta^{(n)}((i,j) \in \op{Inv}(\Sigma_n)) = \frac{1}{(\theta+n-1)(\theta+n-2)} \left(\theta (n-(j-i)) + \binom{n-1}{2} + j-i-1\right).
	\end{align}
	By partial fraction decomposition, this is equal to
	\begin{align*}
		p_{i,j}(\theta):={}&\frac{(n-1)(n-(j-i)) - \binom{n-1}{2} - (j-i-1)}{\theta+n-1} + \frac{\binom{n-1}{2} + j-i-1 - (n-2)(n-(j-i))}{\theta+n-2}\\
		={}& \frac{n(n-2(j-i)+1)}{2 (\theta+n-1)} - \frac{(n-1)(n-2(j-i))}{2 (\theta+n-2)},
	\end{align*}
	proving~\eqref{eqn:P}.
	Setting~$l:=j-i$ and differentiating yields
	$$p_{i,j}'(\theta) = \frac{1}{2}\left(\frac{(n-1)(n-2l)}{(\theta+n-2)^2} - \frac{n(n-2l+1)}{(\theta+n-1)^2}\right).$$ Noting that~$[0,\infty) \to [0,\infty)$, $\theta \mapsto \frac{\theta+n-1}{\theta+n-2}$ is decreasing, we get
	\begin{align*}
		(n-1)(n-2l) \frac{(\theta+n-1)^2}{(\theta+n-2)^2} - n(n-2l+1)
		\le {}&  (n-2l) \frac{(n-1)^3}{(n-2)^2} - n(n-2l+1)\\
		={}& \frac{2l+(2l-5)n - (2l-3)n^2}{(n-2)^2}
		\begin{cases}
			\lne 0, \quad& l \ge 2,\\
			\gne 0, \quad& l = 1,
		\end{cases}
	\end{align*}
	with equality in the first line for~$\theta=0$. Hence~$p_{i,j}(\theta)$ is (strictly) deceasing in~$\theta \ge 0$ if and only if~$l \ge 2$ holds.
	
	$(b)$ Using diagonal summation and well-known summation formulas, we get
	\begin{equation}
		\label{eqn:sum(j-i)}
		\sum_{i,j \in [n],\, i \lne j} (j-i) = \sum_{k=1}^{n-1} (n-k)k = n \frac{n(n-1)}{2} - \frac{(n-1)n(2n-1)}{6} = \binom{n+1}{3}.
	\end{equation}
	Combining this with~\eqref{eqn:Pfirst}, we deduce
	\begin{align*}
		&\E_\theta^{(n)}\abs\op{Inv}(\Sigma_n)\abs\\
		&\qquad= \sum_{i,j \in [n],\, i \lne j} \p_\theta^{(n)} ((i,j) \in \op{Inv}(\Sigma_n))\\
		&\qquad= \frac{1}{(\theta+n-1)(\theta+n-2)} \sum_{i,j \in [n],\, i \lne j} \left(\theta (n-(j-i)) + \binom{n-1}{2} + j-i-1\right)\\
		&\qquad= \frac{1}{(\theta+n-1)(\theta+n-2)} \left( \theta \left(\binom{n}{2} n - \binom{n+1}{3}\right) + \binom{n}{2}\binom{n-1}{2} + \binom{n+1}{3} - \binom{n}{2}\right)\\
		&\qquad= \frac{1}{(\theta+n-1)(\theta+n-2)} \left(\frac{\theta}{4} \binom{2n}{3} + \frac{3n-1}{2} \binom{n}{3}\right).
	\end{align*}
	By partial fraction decomposition, this equals
	\begin{align*}
		g(\theta):=\frac{\frac{n-1}{4} \binom{2n}{3} - \frac{3n-1}{2} \binom{n}{3}}{\theta+n-1} + \frac{\frac{3n-1}{2} \binom{n}{3} - \frac{n-2}{4} \binom{2n}{3}}{\theta+n-2}
		={}& \frac{(n+1)n^2(n-1)}{12 (\theta+n-1)} - \frac{n(n-1)^2(n-2)}{12 (\theta+n-2)},
	\end{align*}
	proving~\eqref{eqn:E}. Differentiation yields
	$$g^{(m)}(\theta) = (-1)^m \frac{m!}{12} n(n-1) \left(\frac{(n+1)n}{(\theta+n-1)^{m+1}} - \frac{(n-1)(n-2)}{(\theta+n-2)^{m+1}}\right), \quad m \in \N.$$ Proceeding as in the proof of part~$(a)$, we get
	\begin{align*}
		(n+1)n - (n-1)(n-2)\frac{(\theta+n-1)^2}{(\theta+n-2)^2}
		\ge{}& (n+1)n - \frac{(n-1)^3}{n-2}
		= \frac{n(2n-5)+1}{n-2}
		\gne 0,
	\end{align*}
	proving~$g' \lne 0$. Similarly, we obtain
	\begin{align*}
		(n+1)n - (n-1)(n-2)\frac{(\theta+n-1)^3}{(\theta+n-2)^3}
		\ge{}& (n+1)n - \frac{(n-1)^4}{(n-2)^2}\\
		={}& \frac{n(n-2)(n-4)-1}{(n-2)^2}
		\begin{cases}
			\gne 0, \quad& n \ge 5,\\
			\lne 0, \quad& n \in \{3,4\},
		\end{cases}
	\end{align*}
	with equality in the first line for~$\theta=0$. Hence~$g(\theta)$ is (strictly) convex in~$\theta \ge 0$ if and only if we have~$n \ge 5$.
\end{proof}

\begin{proof}[Proof of Corollary~\ref{cor:asympExpan}]
	We start by observing
	$$\frac{1}{\theta+n-m} = \frac{1}{n} \cdot \frac{1}{1-\frac{m-\theta}{n}} = \frac{1}{n} \sum_{k=0}^\infty \frac{(m-\theta)^k}{n^k}, \qquad m \in \{1,2\},\ n \gne \theta.$$
	
	$(a)$ Setting~$(l_n)_n:=(j_n-i_n)_n$ and using~\eqref{eqn:P}, we obtain
	\begin{align*}
		&\p_\theta^{(n)}((i_n,j_n) \in \op{Inv}(\Sigma_n))\\
		&\qquad= \frac{n(n-2l_n+1)}{2 (\theta+n-1)} - \frac{n(n-2l_n)}{2 (\theta+n-2)} + \frac{n-2l_n}{2 (\theta+n-2)}\\
		&\qquad= \frac{n-2l_n+1}{2} \sum_{k=0}^\infty \frac{(1-\theta)^k}{n^k} - \frac{n-2l_n}{2} \sum_{k=0}^\infty \frac{(2-\theta)^k}{n^k} + \frac{n-2l_n}{2 n} \sum_{k=0}^\infty \frac{(2-\theta)^k}{n^k}\\
		&\qquad= \frac{1}{2} + \frac{(1-\theta)l_n}{n^2} - \frac{(\theta-1)(\theta-2)}{2n^2} + \frac{(2\theta-3)(\theta-1)l_n}{n^3}
		+ O\left(\frac{1}{n^3}\right),
		\qquad n \to \infty,
	\end{align*}
	where the~$O$-term depends on~$\theta$.
	
	$(b)$ Similarly equation~\eqref{eqn:E} implies
	\begin{align*}
		\E_\theta^{(n)}\abs\op{Inv}(\Sigma_n)\abs
		={}& \frac{(n+1)n^2(n-1)}{12 (\theta+n-1)} - \frac{n(n-1)^2(n-2)}{12 (\theta+n-2)}\\
		={}& \frac{n^3-n}{12} \sum_{k=0}^\infty \frac{(1-\theta)^k}{n^k} - \frac{n^3-4n^2+5n-2}{12} \sum_{k=0}^\infty \frac{(2-\theta)^k}{n^k}\\
		={}& \frac{n(n-1)}{4} + \frac{1-\theta}{6} n + \frac{\theta(\theta-1)}{12} + O\left(\frac{1}{n^2}\right), \qquad n \to \infty,
	\end{align*}
	where the~$O$-term depends on~$\theta$.
\end{proof}


\begin{proof}[Elementary proof of~\eqref{eqn:asmpP} and~\eqref{eqn:asmpE}]
	By the definition of the Ewens sampling distribution, we have~$$\p_\theta^{(n)} ((i,j) \in \op{Inv}(\Sigma_n)) = \frac{1}{\theta^{(n)}} \sum_{k=1}^n a^{i,j}_{n,k} \theta^k,$$ with~$$a^{i,j}_{n,k}:=\abs \{\pi \in S_n: N(\pi)=k,\, (i,j) \in \op{Inv}(\pi)\}\abs, \quad k \in [n],$$ being the number of permutations in~$S_n$ with~$k$ cycles such that~$(i,j)$ is an inversion. Clearly,~$a^{i,j}_{n,n}=0$ holds, showing $$\lim_{\theta \to \infty} \theta \p_\theta^{(n)} ((i,j) \in \op{Inv}(\Sigma_n)) = a^{i,j}_{n,n-1}.$$ Further, we observe that~$\pi \in S_n$ with~$N(\pi)=n-1$ has exactly~$n-2$ fixed points and swaps the remaining two elements of~$[n]$. Thus any~$\pi \in S_n$ with~$N(\pi)=n-1$ satisfies~$(i,j) \in \op{Inv}(\pi)$ if and only if it either swaps~$i$ with some~$k \in \{j+1,\dots,n\}$ or it swaps~$j$ with some~$k \in \{1,\dots,i\}$. We deduce~$a^{i,j}_{n,n-1} = n-j+i$, proving~\eqref{eqn:asmpP}. Regarding~\eqref{eqn:asmpE}, we observe ~$$\E_\theta^{(n)} \abs\op{Inv}(\Sigma_n)\abs = \frac{1}{\theta^{(n)}} \sum_{k=1}^n b_{n,k} \theta^k,$$ with~$b_{n,k}:=\sum_{i,j \in [n],\, i \lne j} a_{n,k}^{i,j}$ for any~$k \in [n]$. In particular,~$b_{n,n}=0$ holds, entailing $$\lim_{\theta \to \infty} \theta \E_\theta^{(n)} \abs\op{Inv}(\Sigma_n)\abs = b_{n,n-1} = \sum_{i,j \in [n],\, i \lne j} (n-j+i) = \frac{1}{4} \binom{2n}{3},$$ where the final equality is derived as in the proof of Theorem~\ref{thm:uncond}$(b)$.
\end{proof}

\section{Proofs of the Results in Subsection~\ref{sec2}}
\label{secP2}
\begin{proof}[Proof of Proposition~\ref{prop:pFix}]
	Let~$n \in \N$ and~$m \in \{0,\dots,n-1\}$.	By symmetry and a basic Chinese restaurant argument, we have
	\begin{align}
		\label{eqn:fixDer}
		\p_\theta^{(n)}(D_{n,m})
		={}& \binom{n}{m} \p_\theta^{(n)}(\op{Fix}(\Sigma_n)=\{n-m+1,\dots,n\}) \notag\\
		={}& \binom{n}{m} \p_\theta^{(n-m)}(D_{n-m,0}) \prod_{l=n-m}^{n-1} \frac{\theta}{\theta+l}
		= \binom{n}{m} \frac{\theta^m\theta^{(n-m)}}{\theta^{(n)}}\p_\theta^{(n-m)}(D_{n-m,0}),
	\end{align}
	proving the first claimed equality.
	For any~$A \subseteq [n]$ with~$k:=\abs A\abs$, a similar approach yields
	\begin{align*}
		\p_\theta^{(n)}(A \subseteq \op{Fix}(\Sigma_n))
		={}& \p_\theta^{(n)}(\{n-k+1,\dots,n\}\subseteq \op{Fix}(\Sigma_n))
		= \prod_{l=n-k}^{n-1} \frac{\theta}{\theta+l}
		= \frac{\theta^k \theta^{(n-k)}}{\theta^{(n)}}.
	\end{align*}
	Together with the inclusion-exclusion principle, we deduce
	\begin{align*}
		\p_\theta^{(n)}(D_{n,0})
		={}& 1- \p_\theta^{(n)}(\op{Fix}(\Sigma_n) \neq \emptyset)\\
		={}& 1- \sum_{k=1}^n (-1)^{k-1} \sum_{A \subseteq [n], \abs A\abs = k} \p_\theta^{(n)}(A \subseteq \op{Fix}(\Sigma_n))
		= \sum_{k=0}^n (-1)^k \binom{n}{k} \frac{\theta^k \theta^{(n-k)}}{\theta^{(n)}}.
	\end{align*}
	It remains to insert this into~\eqref{eqn:fixDer} and simplify the resulting expression to obtain the second claimed equality. The special case~$m=n$ is trivial.
\end{proof}

The main ingredient towards proving the asymptotic results of Subsection~\ref{sec2} is the following proposition. In view of Proposition~\ref{prop:pFix} and equation~\eqref{eqn:pFixTriv}, it provides the exact formulas announced in Subsection~\ref{sec2}:
\begin{prop}
	\label{prop:exact}
	\begin{enumerate}[$(a)$]
		\item We have
		\begin{align}
			\label{eqn:Pexact}
			\p_\theta^{(n)}\big((&i,j) \in \op{Inv}(\Sigma_n)\big\abs D_{n,m}\big)\\
			={}& \frac{1}{(\theta+n-1)(\theta+n-2)} \cdot\notag\\
			&\Bigg[
			\frac{\theta (n-m-1)(n-1-(j-i))}{n-2} \cdot \frac{\p_\theta^{(n-2)}(D_{n-2,m-1})}{\p_\theta^{(n)}(D_{n,m})}\notag\\
			&+ \left(\frac{\theta m(n-1-(j-i))}{n-2} + \theta + (n-m-2)\left(\frac{n-m-3}{2} + \frac{n+j-i-3}{n-2}\right) \right) \frac{\p_\theta^{(n-2)}(D_{n-2,m})}{\p_\theta^{(n)}(D_{n,m})}\notag\\
			&+ (m+1)\left((n-m-3) + \frac{n+j-i-3}{n-2}\right) \frac{\p_\theta^{(n-2)}(D_{n-2,m+1})}{\p_\theta^{(n)}(D_{n,m})}\notag\\
			&+ \frac{(m+2)(m+1)}{2} \cdot \frac{\p_\theta^{(n-2)}(D_{n-2,m+2})}{\p_\theta^{(n)}(D_{n,m})}\Bigg].\notag
		\end{align}
		\item The expected number of inversions under~$\p_\theta^{(n)}(\,\cdot\,\abs D_{n,m})$ is given by
		\begin{align}
			\label{eqn:Eexact}
			\E_\theta^{(n)}\big[\abs\op{Inv}(\Sigma_n)\abs \big\abs D_{n,m}\big]
			={}& \frac{n(n-1)}{(\theta+n-1)(\theta+n-2)} \cdot\\
			&\Bigg[
			\frac{\theta(n-m-1)}{3} \cdot \frac{\p_\theta^{(n-2)}(D_{n-2,m-1})}{\p_\theta^{(n)}(D_{n,m})}\notag\\
			&+ \left(\frac{\theta m}{3} + \frac{\theta}{2} + (n-m-2)\left( \frac{n-m-3}{4} + \frac{2}{3}\right) \right) \frac{\p_\theta^{(n-2)}(D_{n-2,m})}{\p_\theta^{(n)}(D_{n,m})}\notag\\
			&+ (m+1)\left(\frac{n-m-3}{2} + \frac{2}{3}\right) \frac{\p_\theta^{(n-2)}(D_{n-2,m+1})}{\p_\theta^{(n)}(D_{n,m})}\notag\\
			&+ \frac{(m+2)(m+1)}{4} \cdot \frac{\p_\theta^{(n-2)}(D_{n-2,m+2})}{\p_\theta^{(n)}(D_{n,m})}\Bigg].\notag
		\end{align}
	\end{enumerate}
\end{prop}
\begin{proof}$(a)$
	At the top level, we proceed similar to the proof of Theorem~\ref{thm:uncond}. Let~$\Sigma_n$,~$\Sigma_{n-1}$ and~$\Sigma_{n-2}$ be the (consistent) random permutations of~$[n]$, of~$[n] \setminus \{j\}$ and of~$[n] \setminus \{i,j\}$, respectively, arising from a Chinese restaurant construction, modified such that the persons (numbered from~$1$ to~$n$) arrive in the order $$1,\dots,i-1,i+1,\dots,j-1,j+1,\dots,n,i,j.$$ As above, we observe~$\Sigma_n \sim \p_\theta^{(n)}$ under the probability measure encoding the random seating choices of the~$n$ persons, which we denote by~$\p_\theta^{(n)}$ as well. Since we are conditioning on~$D_{n,m}$, we are permanently dealing with the situation that~$\abs\op{Fix}(\Sigma_n)\abs=m$ persons are sitting alone at tables when all~$n$ people are seated. This, in particular, implies~$$\abs\op{Fix}(\Sigma_{n-2})\abs \in \{m-2,\dots,m+2\}.$$ We consider each of the five options separately, thus splitting~$$A:=\{(i,j) \in \op{Inv}(\Sigma_n)\} \cap D_{n,m}=\{\Sigma_n^{-1}(i) \gne \Sigma_n^{-1}(j), \abs\op{Fix}(\Sigma_n)\abs = m\}$$ into five sub-events. We further partition these sub-events into a total of~$13$ sub-sub-events and analyze each of them using the modified Chinese restaurant construction.
	For convenience, we call a table \textit{empty}, a \textit{singleton} or \textit{crowded} according to whether~$0$, $1$ or at least~$2$ persons are seated there.
	
	First consider the situation with~$m-2$ singleton tables when~$n-2$ persons are seated. Then persons~$i$ and~$j$ both would have to choose empty tables, which, however, implies that~$(i,j)$ cannot be an inversion. Hence we have
	\begin{align}
		\label{eqn:m-2}
		&\p_\theta^{(n)}(A \cap \{\abs\op{Fix}(\Sigma_{n-2})\abs = m-2\}=0.
	\end{align}
	
	Now consider the situation with~$m-1$ singleton tables when~$n-2$ persons are seated. We make a case distinction with respect to the seating choice of~$i$. We start by assuming that~$i$ chooses an empty table. Then~$j$ must choose a crowded table to ensure that we end up with~$m$ singleton tables. By symmetry, the~$m-1$ persons sitting at singleton tables when~$i$ arrives (i.e., the fixed points of~$\Sigma_{n-2}$) are distributed with equal probability among the~$n-2$ seated persons. Using this symmetry property in the fourth equality, we get
	\begin{align*}
		&\p_\theta^{(n)}(A \cap \{\abs\op{Fix}(\Sigma_{n-2})\abs = m-1, \Sigma_{n-1}(i)=i\})\\
		&\qquad= \p_\theta^{(n)}(\Sigma_n^{-1}(j) \not\in \op{Fix}(\Sigma_{n-2}), i \gne \Sigma_n^{-1}(j), \Sigma_{n-1}(i)=i, \abs\op{Fix}(\Sigma_{n-2})\abs = m-1)\\
		&\qquad= \sum_{k=1}^{i-1} \p_\theta^{(n)}(\Sigma_n^{-1}(j)=k, \Sigma_{n-1}(i)=i, k \not\in \op{Fix}(\Sigma_{n-2}), \abs\op{Fix}(\Sigma_{n-2})\abs = m-1)\\
		&\qquad= \sum_{k=1}^{i-1} \frac{1}{\theta+n-1} \cdot \frac{\theta}{\theta+n-2} \p_\theta^{(n-2)}(\{k \not\in \op{Fix}(\Sigma_{n-2})\} \cap D_{n-2,m-1})\\
		&\qquad= \sum_{k=1}^{i-1} \frac{1}{\theta+n-1} \cdot \frac{\theta}{\theta+n-2} \cdot \frac{(n-2)-(m-1)}{n-2} \p_\theta^{(n-2)}(D_{n-2,m-1})\\
		&\qquad= \frac{i-1}{\theta+n-1}\cdot \frac{n-m-1}{n-2} \cdot \frac{\theta}{\theta+n-2} \p_\theta^{(n-2)}(D_{n-2,m-1}).
	\end{align*}
	Throughout the rest of the proof, we make extensive use of symmetry arguments similar to the last four lines without repeating the details. Let us now assume that~$i$ does not choose an empty table. Then necessarily~$i$ must choose a crowded table while~$j$ must choose an empty table. Noting that this implies~$\Sigma_n^{-1}(i)= \Sigma_{n-1}^{-1}(i)$, a similar calculation yields
	\begin{align*}
		&\p_\theta^{(n)}(A \cap \{\abs\op{Fix}(\Sigma_{n-2})\abs = m-1, \Sigma_{n-1}(i)\neq i\})\\
		&\qquad= \p_\theta^{(n)}(\Sigma_n(j)=j, \Sigma_{n-1}^{-1}(i) \not \in \op{Fix}(\Sigma_{n-2}), \Sigma_{n-1}^{-1}(i) \gne j, \abs\op{Fix}(\Sigma_{n-2})\abs = m-1)\\
		&\qquad= \frac{\theta}{\theta+n-1} \cdot \frac{n-m-1}{n-2}\cdot \frac{n-j}{\theta+n-2} \p_\theta^{(n-2)}(D_{n-2,m-1}).
	\end{align*}
	We deduce
	\begin{align}
		\label{eqn:m-1}
		\p_\theta^{(n)}(A \cap \{\abs\op{Fix}(\Sigma_{n-2})\abs = m-1\})
		= \frac{\theta(n-m-1)(n-1-(j-i))}{(\theta+n-1)(\theta+n-2)(n-2)} \p_\theta^{(n-2)}(D_{n-2,m-1}).
	\end{align}
	
	Next consider the situation with~$m$ singleton tables when~$n-2$ persons are seated. We make a case distinction with respect to the seating choices of~$i$ and~$j$. First, assume that~$i$ chooses an empty table and that~$j$ joins~$i$. Then~$(i,j)$ is an inversion and we have
	\begin{align*}
		&\p_\theta^{(n)}(A \cap \{\abs\op{Fix}(\Sigma_{n-2})\abs = m, \Sigma_{n-1}^{-1}(i) = i, \Sigma_n^{-1}(j)=i\})\\
		&\qquad= \p_\theta^{(n)}(\Sigma_n^{-1}(j)=i, \Sigma_{n-1}^{-1}(i)=i, \abs\op{Fix}(\Sigma_{n-2})\abs = m)\\
		&\qquad= \frac{1}{\theta+n-1}\cdot \frac{\theta}{\theta+n-2} \p_\theta^{(n-2)}(D_{n-2,m}).
	\end{align*}
	Second, assume that~$i$ chooses an empty table and that~$j$ does not join~$i$. Then~$j$ must sit down at one of the other singleton tables. Noting~$\Sigma_n^{-1}(i)=i$ and using our usual symmetry argument, we get
	\begin{align*}
		&\p_\theta^{(n)}(A \cap \{\abs\op{Fix}(\Sigma_{n-2})\abs = m, \Sigma_{n-1}^{-1}(i) = i, \Sigma_n^{-1}(j)\neq i\})\\
		&\qquad= \p_\theta^{(n)}(\Sigma_n^{-1}(j) \in \op{Fix}(\Sigma_{n-2}), i \gne \Sigma_n^{-1}(j), \Sigma_{n-1}^{-1}(i)=i, \abs\op{Fix}(\Sigma_{n-2})\abs = m)\\
		&\qquad= \frac{m}{n-2} \cdot \frac{i-1}{\theta+n-1}\cdot \frac{\theta}{\theta+n-2} \p_\theta^{(n-2)}(D_{n-2,m}).
	\end{align*}
	Third, assume that~$i$ chooses a singleton table. Then~$j$ must choose an empty table and we obtain
	\begin{align*}
		&\p_\theta^{(n)}(A \cap \{\abs\op{Fix}(\Sigma_{n-2})\abs = m, \Sigma_{n-1}^{-1}(i) \in \op{Fix}(\Sigma_{n-2})\})\\
		&\qquad= \p_\theta^{(n)}(\Sigma_n(j)=j, \Sigma_{n-1}^{-1}(i) \in \op{Fix}(\Sigma_{n-2}), \Sigma_{n-1}^{-1}(i) \gne j,  \abs\op{Fix}(\Sigma_{n-2})\abs = m)\\
		&\qquad= \frac{\theta}{\theta+n-1}\cdot \frac{m}{n-2} \cdot \frac{n-j}{\theta+n-2} \p_\theta^{(n-2)}(D_{n-2,m}).
	\end{align*}
	Fourth, assume that~$i$ chooses a crowded table and that~$j$ sits down to the right of~$i$. Then the original right neighbor of~$i$ becomes the right neighbor of~$j$ (i.e.,~$\Sigma_{n-1}^{-1}(i)= \Sigma_n^{-1}(j)$). We obtain
	\begin{align*}
		&\p_\theta^{(n)}(A \cap \{\abs\op{Fix}(\Sigma_{n-2})\abs = m, \Sigma_{n-1}^{-1}(i) \not\in \op{Fix}(\Sigma_{n-2}) \cup \{i\}, \Sigma_n(j)=i\})\\
		&\qquad= \p_\theta^{(n)}(\Sigma_n(j)=i, \Sigma_{n-1}^{-1}(i) \not\in \op{Fix}(\Sigma_{n-2}), j \gne \Sigma_{n-1}^{-1}(i) \neq i, \abs\op{Fix}(\Sigma_{n-2})\abs = m)\\
		&\qquad= \frac{1}{\theta+n-1}\cdot \frac{n-2-m}{n-2}\cdot\frac{j-2}{\theta+n-2} \p_\theta^{(n-2)}(D_{n-2,m}).
	\end{align*}
	Fifth, assume that~$i$ chooses a crowded table and that~$j$ sits down to the left of~$i$. Then the right neighbor of~$i$ remains the same after~$j$ is seated (i.e.,~$\Sigma_{n-1}^{-1}(i)= \Sigma_n^{-1}(i)$). We obtain
	\begin{align*}
		&\p_\theta^{(n)}(A \cap \{\abs\op{Fix}(\Sigma_{n-2})\abs = m, \Sigma_{n-1}^{-1}(i) \not\in \op{Fix}(\Sigma_{n-2}) \cup \{i\}, \Sigma_n^{-1}(j)=i\})\\
		&\qquad= \p_\theta^{(n)}(\Sigma_n^{-1}(j)=i, \Sigma_{n-1}^{-1}(i) \not\in \op{Fix}(\Sigma_{n-2}), \Sigma_{n-1}^{-1}(i) \gne i, \abs\op{Fix}(\Sigma_{n-2})\abs = m)\\
		&\qquad= \frac{1}{\theta+n-1}\cdot \frac{n-2-m}{n-2}\cdot\frac{n-i-1}{\theta+n-2} \p_\theta^{(n-2)}(D_{n-2,m}).
	\end{align*}
	Sixth and lastly, assume that~$i$ chooses a crowded table and that~$j$ does not sit down next to~$i$. Then~$j$ must have joined a crowded table as well.
	Further,~$i$ must have chosen to sit to the left of some~$k_i \in [n] \setminus \{i,j\}$ (i.e.,~$\Sigma_{n-1}^{-1}(i)=k_i$). This~$k_i$ remains the right neighbor~$i$ (i.e.,~$\Sigma_n^{-1}(i)=k_i$). Hence~$(i,j)$ is an inversion in this scenario if and only if~$j$ sits down to the left of a person~$k_j \in [n] \setminus \{i,j\}$ with~$k_i \gne k_j$.
	For such~$k_i,k_j \in [n] \setminus \{i,j\}$ with~$k_i \gne k_j$, a (slightly more elaborate) symmetry argument yields
	\begin{align*}
		&\p_\theta^{(n)}(\abs\op{Fix}(\Sigma_n)\abs =m ,\abs\op{Fix}(\Sigma_{n-2})\abs = m, \Sigma_n^{-1}(j) = k_j, \Sigma_{n-1}^{-1}(i)=k_i)\\
		&\qquad= \p_\theta^{(n)}(\Sigma_n^{-1}(j) = k_j, \Sigma_{n-1}^{-1}(i)=k_i, k_i \not\in \op{Fix}(\Sigma_{n-2}), k_j \not\in \op{Fix}(\Sigma_{n-2}), \abs\op{Fix}(\Sigma_{n-2})\abs = m)\\
		&\qquad= \frac{1}{\theta+n-1}\cdot \frac{1}{\theta+n-2} \cdot \frac{(n-2-m)(n-3-m)}{(n-2)(n-3)} \1_{n \gne 3} \p_\theta^{(n-2)}(D_{n-2,m+1}).
	\end{align*}
	Noting
	\begin{equation}
		\label{eqn:triangle}
		\sum_{k_i,k_j \in [n] \setminus \{i,j\}} \1_{k_i \gne k_j} = \frac{(n-2)(n-3)}{2},
	\end{equation}
	we deduce
	\begin{align*}
		&\p_\theta^{(n)}(A \cap \{\abs\op{Fix}(\Sigma_{n-2})\abs = m, \Sigma_n(j) \neq i, \Sigma_n^{-1}(j) \neq i)\})\\
		&\qquad= \sum_{k_i,k_j \in [n] \setminus \{i,j\}} \1_{k_i \gne k_j} \p_\theta^{(n)}(\abs\op{Fix}(\Sigma_n)\abs =m ,\abs\op{Fix}(\Sigma_{n-2})\abs = m, \Sigma_n^{-1}(j) = k_j, \Sigma_{n-1}^{-1}(i)=k_i)\\
		&\qquad= \frac{(n-m-2)(n-m-3)}{2(\theta+n-1)(\theta+n-2)} \p_\theta^{(n-2)}(D_{n-2,m}).
	\end{align*}
	(The indicator~$\1_{n \gne 3}$ can be omitted since~$n=3$ implies~$(n-m-2)(n-m-3)=0$.)
	Combining this with the probabilities covering the other five possible seating arrangements treated above, we get
	\begin{align}
		\label{eqn:m}
		&\p_\theta^{(n)}(A \cap \{\abs\op{Fix}(\Sigma_{n-2})\abs = m\})\\
		={}& \frac{\p_\theta^{(n-2)}(D_{n-2,m})}{(\theta+n-1)(\theta+n-2)} \left(\frac{\theta m(n-1-(j-i))}{n-2} + \theta + (n-2-m)\left(\frac{n+j-i-3}{n-2} + \frac{n-m-3}{2}\right) \right).\notag
	\end{align}
	
	Next consider the situation with~$m+1$ singleton tables when~$n-2$ persons are seated. We make a case distinction with respect to the seating choice of~$j$. First, assume that~$j$ chooses to sit to the right of~$i$. Then~$i$ must have joined a singleton table and the original seatmate of~$i$ becomes the right neighbor of~$j$ (i.e.,~$\Sigma_{n-1}^{-1}(i)= \Sigma_n^{-1}(j)$). Noting that~$i$ joining a singleton table implies~$\Sigma_{n-1}^{-1}(i)\neq i$ and using our usual symmetry argument, we get
	\begin{align*}
		&\p_\theta^{(n)}(A \cap \{\abs\op{Fix}(\Sigma_{n-2})\abs = m+1, \Sigma_n(j) = i\})\\
		&\qquad= \p_\theta^{(n)}(\Sigma_n(j)=i, \Sigma_{n-1}^{-1}(i) \in \op{Fix}(\Sigma_{n-2}), j \gne \Sigma_{n-1}^{-1}(i) \neq i, \abs\op{Fix}(\Sigma_{n-2})\abs = m+1)\\
		&\qquad= \frac{1}{\theta+n-1}\cdot \frac{m+1}{n-2}\cdot \frac{j-2}{\theta+n-2}\p_\theta^{(n-2)}(D_{n-2,m+1}).
	\end{align*}
	Second, assume that~$j$ chooses to sit to the left of~$i$. Then, again,~$i$ must have joined a singleton table and the original seatmate of~$i$ remains its right neighbor (i.e.,~$\Sigma_{n-1}^{-1}(i)= \Sigma_n^{-1}(i)$). We obtain
	\begin{align*}
		&\p_\theta^{(n)}(A \cap \{\abs\op{Fix}(\Sigma_{n-2})\abs = m+1, \Sigma_n^{-1}(j) = i\})\\
		&\qquad= \p_\theta^{(n)}(\Sigma_n^{-1}(j)=i, \Sigma_{n-1}^{-1}(i) \in \op{Fix}(\Sigma_{n-2}), \Sigma_{n-1}^{-1}(i) \gne i, \abs\op{Fix}(\Sigma_{n-2})\abs = m+1)\\
		&\qquad= \frac{1}{\theta+n-1}\cdot \frac{m+1}{n-2}\cdot \frac{n-i-1}{\theta+n-2}\p_\theta^{(n-2)}(D_{n-2,m+1}).
	\end{align*}
	Third and lastly, assume that~$j$ chooses not to sit next to~$i$. We observe that~$j$ cannot choose an empty table. Not being allowed to sit at an empty table either,~$i$ must have chosen to sit to the left of some~$k_i \in [n] \setminus \{i,j\}$ (i.e.,~$\Sigma_{n-1}^{-1}(i)=k_i$).
	This~$k_i$ remains the right neighbor~$i$ (i.e.,~$\Sigma_n^{-1}(i)=k_i$). Hence~$(i,j)$ is an inversion in this scenario if and only if~$j$ sits down to the left of a person~$k_j \in [n] \setminus \{i,j\}$ with~$k_i \gne k_j$.
	For such~$k_i,k_j \in [n] \setminus \{i,j\}$, a symmetry argument yields
	\begin{align*}
		&\p_\theta^{(n)}(\abs\op{Fix}(\Sigma_n)\abs =m ,\abs\op{Fix}(\Sigma_{n-2})\abs = m+1, \Sigma_n^{-1}(j) = k_j, \Sigma_{n-1}^{-1}(i)=k_i)\\
		&\qquad= \p_\theta^{(n)}(\Sigma_n^{-1}(j) = k_j, \Sigma_{n-1}^{-1}(i)=k_i, k_i \in \op{Fix}(\Sigma_{n-2}), k_j \not\in \op{Fix}(\Sigma_{n-2}), \abs\op{Fix}(\Sigma_{n-2})\abs = m+1)\\
		&\qquad\quad+ \p_\theta^{(n)}(\Sigma_n^{-1}(j) = k_j, \Sigma_{n-1}^{-1}(i)=k_i, k_i \not\in \op{Fix}(\Sigma_{n-2}), k_j \in \op{Fix}(\Sigma_{n-2}), \abs\op{Fix}(\Sigma_{n-2})\abs = m+1)\\
		&\qquad= 2\cdot \frac{1}{\theta+n-1}\cdot \frac{1}{\theta+n-2} \cdot \frac{(m+1)((n-2)-(m+1))}{(n-2)(n-3)} \1_{n \gne 3} \p_\theta^{(n-2)}(D_{n-2,m+1}).
	\end{align*}
	Using~\eqref{eqn:triangle},	we deduce
	\begin{align*}
		&\p_\theta^{(n)}(A \cap \{\abs\op{Fix}(\Sigma_{n-2})\abs = m+1, \Sigma_n(j) \neq i, \Sigma_n^{-1}(j) \neq i)\})\\
		&\qquad= \sum_{k_i,k_j \in [n] \setminus \{i,j\}} \1_{k_i \gne k_j} \p_\theta^{(n)}(\abs\op{Fix}(\Sigma_n)\abs =m ,\abs\op{Fix}(\Sigma_{n-2})\abs = m+1, \Sigma_n^{-1}(j) = k_j, \Sigma_{n-1}^{-1}(i)=k_i)\\
		&\qquad= \frac{(m+1)(n-m-3)}{(\theta+n-1)(\theta+n-2)} \p_\theta^{(n-2)}(D_{n-2,m+1}).
	\end{align*}
	(The indicator~$\1_{n \gne 3}$ can be omitted since~$n=3$ implies~$(n-m-3)\p_\theta^{(n-2)}(D_{n-2,m+1})=0$.)
	Combining this with the probabilities covering the other two seating choices of~$j$ treated above, we get
	\begin{align}
		\label{eqn:m+1}
		&\p_\theta^{(n)}(A \cap \{\abs\op{Fix}(\Sigma_{n-2})\abs = m+1\})\\
		&\qquad= \frac{m+1}{(\theta+n-1)(\theta+n-2)}\left(\frac{n+j-i-3}{n-2}+ (n-m-3)\right) \p_\theta^{(n-2)}(D_{n-2,m+1}).\notag
	\end{align}
	
	Finally, consider the situation with~$m+2$ singleton tables when~$n-2$ persons are seated. Then persons~$i$ and~$j$ have to join two distinct such singleton tables. By symmetry, this leads to~$(i,j)$ being an inversion with conditional probability~$\frac{1}{2}$. Thus we have
	\begin{align}
		\label{eqn:m+2}
		&\p_\theta^{(n)}(A \cap \{\abs\op{Fix}(\Sigma_{n-2})\abs = m+2\})
		= \frac{1}{2} \cdot \frac{m+1}{\theta+n-1} \cdot \frac{m+2}{\theta+n-2} \p_\theta^{(n-2)}(D_{n-2,m+2}).
	\end{align}
	
	Dividing equations~\eqref{eqn:m-2},~\eqref{eqn:m-1},~\eqref{eqn:m},~\eqref{eqn:m+1} and~\eqref{eqn:m+2} by~$\p_\theta^{(n)}(D_{n,m})$ and adding them up, we obtain~\eqref{eqn:Pexact}.
	
	$(b)$ We start by observing
	\begin{align*}
		\E_\theta^{(n)}\big[\abs\op{Inv}(\Sigma_n)\abs \big\abs D_{n,m}\big]
		={}& \sum_{i,j \in [n],\, i \lne j} \p_\theta^{(n)}\big((i,j) \in \op{Inv}(\Sigma_n)\big\abs D_{n,m}\big).
	\end{align*}
	Inserting~\eqref{eqn:Pexact} into this equation, using~$\sum_{i,j \in [n],\, i \lne j} (j-i) = \binom{n+1}{3}$ (cf.~\eqref{eqn:sum(j-i)}) and simplifying the result, we end up with~\eqref{eqn:Eexact}.
\end{proof}

\begin{proof}[Proof of Theorem~\ref{thm:condN}]
	Unless explicitly indicated otherwise, all asymptotic expansions in this proof are expansions as~$n \to \infty$. Let~$b \in \N_0$. We note that~$\theta^{(b)}=\frac{\Gamma(\theta+b)}{\Gamma(\theta)}$, where~$\Gamma$ denotes the Gamma function. In view of Proposition~\ref{prop:exact} and Proposition~\ref{prop:pFix}, we start by analyzing the asymptotic behavior of~
	\begin{align}
		\label{eqn:Gdfn}
		G_b(n):=\sum_{k=0}^{n-b} \frac{(-\theta)^k \theta^{(n-b-k)}}{k!(n-b-k)!} = \frac{1}{\Gamma(\theta)} \sum_{k=0}^{n-b} \frac{(-\theta)^k \Gamma(\theta+n-b-k)}{k!\Gamma(1+n-b-k)}.
	\end{align}
	Let~$l_n:=\lfloor \log(n)\rfloor$. Noting that~$\frac{1}{(l_n+1)!}$ decays superpolynomially, we get
	\begin{align}
		\label{eqn:remainder}
		\left\abs \sum_{k=l_n+1}^{n-b} \frac{(-\theta)^k \theta^{(n-b-k)}}{k!(n-b-k)!}\right\abs
		\le{}& \sum_{k=l_n+1}^{n-b} \frac{\theta^k (\lceil\theta\rceil+n-b-k)!}{k!(n-b-k)!}\\
		\le{}& (\lceil \theta\rceil +n)^{\lceil \theta\rceil} \sum_{k=l_n+1}^\infty \frac{\theta^k}{k!}
		\le (\lceil \theta\rceil +n)^{\lceil \theta\rceil} \frac{\e^\theta}{(l_n+1)!}
		= O\left(\frac{1}{n^4}\right).\notag
	\end{align}
	
	We now focus on the truncated sum up to~$k=l_n$.
	For all~$\alpha,\beta \in \mathbb{R}$ and as~$z \to \infty$,~\cite{TE51} yields
	\begin{align*}
		\frac{\Gamma(\alpha+z)}{\Gamma(\beta+z)} = z^{\alpha-\beta} \left[1+ \frac{(\alpha-\beta)(\alpha+\beta-1)}{2z} + \binom{\alpha-\beta}{2} \frac{3(\alpha+\beta-1)^2 -\alpha+\beta-1}{12 z^2} + O\left(\frac{1}{z^3}\right)\right].
	\end{align*}
	This implies
	\begin{align*}
		\frac{\Gamma(\theta+n-b-k)}{\Gamma(1+n-b-k)} = (n-b-k)^{\theta-1} \left[1+ \frac{(\theta-1)\theta}{2(n-b-k)} + \binom{\theta-1}{2} \frac{3 \theta^2-\theta}{12(n-b-k)^2} + O\left(\frac{1}{n^3}\right)\right],
	\end{align*}
	uniformly over~$k \in \{0,\dots,l_n\}$, entailing
	\begin{align*}
		H_b(n):={}&\sum_{k=0}^{l_n} \frac{(-\theta)^k \Gamma(\theta+n-b-k)}{k!\Gamma(1+n-b-k)}\\
		={}& \sum_{k=0}^{l_n} \frac{(-\theta)^k}{k!} (n-b-k)^{\theta-1}\\
		&+ \sum_{k=0}^{l_n} \frac{(-\theta)^k}{k!} \cdot \frac{\theta(\theta-1)}{2} (n-b-k)^{\theta-2}\\
		&+ \sum_{k=0}^{l_n} \frac{(-\theta)^k}{k!} \cdot \frac{(\theta-1)(\theta-2)(3 \theta^2-\theta)}{24}(n-b-k)^{\theta-3}\\
		&+ \sum_{k=0}^{l_n} \frac{(-\theta)^k}{k!} O(n^{\theta-4}).
	\end{align*}
	Observing
	\begin{align*}
		(n-b-k)^{\theta-a}
		={}& n^{\theta-a} \left(1-\frac{b+k}{n}\right)^{\theta-a}\\
		={}& n^{\theta-a} \left[1- \frac{(\theta-a)(b+k)}{n} + \frac{(\theta-a)(\theta-a-1)(b+k)^2}{2n^2} + O\left(\frac{1}{n^3}\right)\right],
	\end{align*}
	uniformly over~$k \in \{0,\dots,l_n\}$ and~$a \in \{1,2,3\}$, we deduce
	\begin{align}
		\label{eqn:H}
		H_b(n)
		={}& n^{\theta-1} \sum_{k=0}^{l_n} \frac{(-\theta)^k}{k!}
		+ n^{\theta-2} \sum_{k=0}^{l_n} \frac{(-\theta)^k}{k!} \left[\frac{\theta(\theta-1)}{2}-(\theta-1)(b+k)\right]\\
		&+ n^{\theta-3} \sum_{k=0}^{l_n} \frac{(-\theta)^k}{k!} \cdot \frac{(\theta-1)(\theta-2)}{2} \left[\frac{3 \theta^2-\theta}{12} - \theta(b+k)
		+ (b+k)^2\right]
		+ O(n^{\theta-4}).\notag
	\end{align}
	Similar to~\eqref{eqn:remainder}, the superpolynomial decay of~$\frac{1}{(l_n+1)!}$ yields
	\begin{align}
		\label{eqn:exp}
		\sum_{k=0}^{l_n} \frac{(-\theta)^k}{k!} ={}& \e^{-\theta} + O\left(\frac{1}{n^3}\right),\notag\\
		\sum_{k=0}^{l_n} \frac{(-\theta)^k}{k!} k ={}& -\theta \sum_{k=0}^{l_n-1} \frac{(-\theta)^k}{k!}= -\theta \e^{-\theta} + O\left(\frac{1}{n^3}\right),\\
		\sum_{k=0}^{l_n} \frac{(-\theta)^k}{k!} k^2 ={}& \theta^2 \sum_{k=0}^{l_n-2} \frac{(-\theta)^k}{k!} - \theta \sum_{k=0}^{l_n-1} \frac{(-\theta)^k}{k!} = \theta(\theta-1) \e^{-\theta} + O\left(\frac{1}{n^3}\right).\notag	
	\end{align}
	Hence~\eqref{eqn:Gdfn}-\eqref{eqn:exp} imply
	\begin{align}
		\label{eqn:G2}
		G_b(n)
		={}& \frac{H_b(n)}{\Gamma(\theta)} + O\left(\frac{1}{n^4}\right)
		= \frac{n^{\theta-1}}{\Gamma(\theta)\e^{\theta}} \bigg[1 + \frac{1}{n} A_1(b) + \frac{1}{n^2} A_2(b) + O\left(\frac{1}{n^3}\right)\bigg],
	\end{align}
	with
	\begin{align}
		\label{eqn:A1}
		A_1(b)
		:={}& \frac{\theta(\theta-1)}{2} - (\theta-1)(b -\theta) = \frac{\theta-1}{2} (3\theta-2b)
	\end{align}
	and
	\begin{align}
		\label{eqn:A2}
		A_2(b)
		:={}& \frac{(\theta-1)(\theta-2)}{2} \left(\frac{3\theta^2-\theta}{12} - \theta (b-\theta) + b^2-2b\theta+\theta(\theta-1)\right)\\
		={}& \frac{(\theta-1)(\theta-2)}{2} \left(b^2 - 3b\theta + \theta \frac{27\theta-13}{12}\right).\notag
	\end{align}
	
	Now let~$a \in \{-1,0,1,2\}$. Using a geometric series expansion to rewrite the denominator, it follows from~\eqref{eqn:G2}, \eqref{eqn:A1} and~\eqref{eqn:A2} that
	\begin{align}
		\label{eqn:Gquot}
		\frac{G_{m+2+a}(n)}{G_m(n)}
		={}& \frac{1 + \frac{1}{n} A_1(m+2+a) + \frac{1}{n^2} A_2(m+2+a) + O\left(\frac{1}{n^3}\right)}{1 + \frac{1}{n} A_1(m) + \frac{1}{n^2} A_2(m) + O\left(\frac{1}{n^3}\right)}\\
		={}& \left[1 + \frac{A_1(m+2+a)}{n} + \frac{A_2(m+2+a)}{n^2} + O\left(\frac{1}{n^3}\right)\right]\notag\\
		&\cdot \left[1 - \frac{A_1(m)}{n} - \frac{A_2(m)}{n^2} + \frac{(A_1(m))^2}{n^2} + O\left(\frac{1}{n^3}\right) \right]\notag\\
		={}& 1+ \frac{B_1(m,a)}{n} + \frac{B_2(m,a)}{n^2} + O\left(\frac{1}{n^3}\right),\notag
	\end{align}
	with
	\begin{align}
		\label{eqn:B1}
		B_1(m,a)
		:={}& A_1(m+2+a) - A_1(m)
		= -(2+a)(\theta-1)
	\end{align}
	and
	\begin{align}
		\label{eqn:B2}
		B_2(m,a)
		:={}& A_2(m+2+a) - A_1(m+2+a)A_1(m) - A_2(m) + (A_1(m))^2\\
		={}& \frac{(2+a)(\theta-1)}{2} ((5+a)\theta - 2(m+2+a)).\notag
	\end{align}
	Combining~\eqref{eqn:Gdfn} and~\eqref{eqn:Gquot} with Proposition~\ref{prop:pFix}, we obtain
	\begin{align}
		\label{eqn:pQuot}
		&\frac{1}{(\theta+n-1)(\theta+n-2)} \cdot \frac{\p_\theta^{(n-2)}(D_{n-2,m+a})}{\p_\theta^{(n)}(D_{n,m})}\\
		&\qquad= \frac{m!\theta^a}{(m+a)!n(n-1)} \cdot\frac{G_{m+2+a}(n)}{G_m(n)}\notag\\
		&\qquad= \frac{m!\theta^a}{(m+a)!n(n-1)} \left[1+ \frac{B_1(m,a)}{n} + \frac{B_2(m,a)}{n^2} + O\left(\frac{1}{n^3}\right)\right].\notag
	\end{align}
	(In the special case~$(m,a)=(0,-1)$, this trivially holds under the convention~$\frac{0!}{(-1)!}:=0$.) To complete the proof of part~$(b)$, it remains to insert~\eqref{eqn:pQuot} into~\eqref{eqn:Eexact}, use~\eqref{eqn:B1} and~\eqref{eqn:B2}, and simplify the resulting expression.
	
	Regarding part~$(a)$, we additionally observe~$$\frac{1}{n-1} = \frac{1}{n}+\frac{1}{n^2}+\frac{1}{n^3}+O\left(\frac{1}{n^4}\right).$$ Thus, using~\eqref{eqn:B1} and~\eqref{eqn:B2}, equation~\eqref{eqn:pQuot} can be rewritten as
	\begin{align*}
		&\frac{1}{(\theta+n-1)(\theta+n-2)} \cdot \frac{\p_\theta^{(n-2)}(D_{n-2,m+a})}{\p_\theta^{(n)}(D_{n,m})}\\
		&\qquad= \frac{m!\theta^a}{(m+a)!} \left[\frac{1}{n^2}+ \frac{1-(2+a)(\theta-1)}{n^3}+ \frac{2+(2+a)(\theta-1)((5+a)\theta-2(m+3+a))}{2n^4}
		+ O\left(\frac{1}{n^5}\right)\right].
	\end{align*}
	Inserting this (the $O(1/n^4)$-term is only needed for~$a=0$) together with~$$\frac{1}{n-2} = \frac{1}{n} + \frac{2}{n^2} + O \left(\frac{1}{n^3}\right)$$ into~\eqref{eqn:Pexact}, it remains to simplify the resulting expression.
\end{proof}

\begin{proof}[Proof of Theorem~\ref{thm:condThetaInf}]
	Let~$l \in \N$.	By the definition of the Ewens sampling distribution, we have~$$\p_\theta^{(l)} (D_{l,0}) = \frac{1}{\theta^{(l)}} \sum_{k=1}^l a_{l,k} \theta^k,$$ with~$$a_{l,k}:=\abs \{\pi \in S_l: N(\pi)=k,\, \op{Fix}(\pi)=\emptyset\}\abs, \quad k \in [l],$$ being the number of derangements in~$S_l$ with exactly~$k$ cycles. Observing~$a_{l,k}=0$ for all~$k \in [l]$ with~$k \gne \frac{l}{2}$, we get
	\begin{equation}
		\label{eqn:a}
			\lim_{\theta \to \infty} \theta^{\lceil l/2\rceil} \p_\theta^{(l)}(D_{l,0}) = a_{l,\lfloor l/2\rfloor}.		
	\end{equation}
	If~$l$ is even, $\pi \in S_l$ is a derangement with~$\lfloor l/2\rfloor$ cycles if and only if it consists of exactly~$l/2$ two-cycles, entailing~$a_{l,\lfloor l/2\rfloor} = (l-1)!!$. If~$l\ge 3$ is odd, $\pi \in S_l$ is a derangement with~$\lfloor l/2\rfloor$ cycles if and only if it consists of exactly~$1$~three-cycle and~$(l-3)/2$ two-cycles. Thus, for any odd~$l$, elementary combinatorics implies~$a_{l,\lfloor l/2\rfloor} = 2\binom{l}{3} (l-4)!!$, with the conventions~$(-1)!!:=1$ and~$\binom{1}{3} (-3)!!:=0$.
	
	For the time being, assume that~$n-m$ is even. Using the first equality in Proposition~\ref{prop:pFix} along with~\eqref{eqn:a} and the subsequent formulas, we obtain
	\begin{align*}
		\lim_{\theta \to \infty} \frac{\p_\theta^{(n-2)}(D_{n-2,m})}{\theta \p_\theta^{(n)}(D_{n,m})}
		={}& \lim_{\theta \to \infty} \frac{\binom{n-2}{m} \theta^{\lceil (n-m-2)/2\rceil} \p_\theta^{(n-m-2)}(D_{n-m-2,0})}{\binom{n}{m} \theta^{\lceil (n-m)/2\rceil} \p_\theta^{(n-m)}(D_{n-m,0})}\\
		={}& \frac{\binom{n-2}{m} (n-m-3)!!}{\binom{n}{m} (n-m-1)!!}
		= \frac{n-m}{n(n-1)}.
	\end{align*}
	(For~$m=n-2$, this is true on account of our conventions.) Similarly (or trivially for~$m=n-2$), we get
	\begin{align*}
		\lim_{\theta \to \infty} \frac{\p_\theta^{(n-2)}(D_{n-2,m+2})}{\theta^2 \p_\theta^{(n)}(D_{n,m})}
		= \frac{(n-m)(n-m-2)}{n(n-1)(m+2)(m+1)}
	\end{align*}
	as well as
	\begin{align*}
		0=\lim_{\theta \to \infty} \frac{\p_\theta^{(n-2)}(D_{n-2,m-1})}{\theta \p_\theta^{(n)}(D_{n,m})}
		= \lim_{\theta \to \infty} \frac{\p_\theta^{(n-2)}(D_{n-2,m+1})}{\theta^2 \p_\theta^{(n)}(D_{n,m})}.
	\end{align*}
	It remains to combine these four limits with Proposition~\ref{prop:exact} and simplify the resulting expressions.
	
	Now assume that~$n-m$ is odd. Then a similar approach yields
	\begin{align*}
		\lim_{\theta \to \infty} \frac{\p_\theta^{(n-2)}(D_{n-2,m-1})}{\theta \p_\theta^{(n)}(D_{n,m})}
		={}& \frac{\binom{n-2}{m-1} \theta^{\lceil (n-m-1)/2\rceil} \p_\theta^{(n-m-1)}(D_{n-m-1,0})}{\binom{n}{m} \theta^{\lceil (n-m)/2\rceil} \p_\theta^{(n-m)}(D_{n-m,0})}\\
		={}& \frac{\binom{n-2}{m-1} (n-m-2)!!}{\binom{n}{m} 2 \binom{n-m}{3} (n-m-4)!!}
		= \frac{3m}{n(n-1)(n-m-1)}
	\end{align*}
	(which is trivial for~$m=0$) and
	\begin{align*}
		\lim_{\theta \to \infty} \frac{\p_\theta^{(n-2)}(D_{n-2,m+1})}{\theta^2 \p_\theta^{(n)}(D_{n,m})}
		= \frac{3}{n(n-1)(m+1)}
	\end{align*}
	as well as
	\begin{align*}
		\lim_{\theta \to \infty} \frac{\p_\theta^{(n-2)}(D_{n-2,m})}{\theta \p_\theta^{(n)}(D_{n,m})}
		={}& \frac{\binom{n-2}{m} \theta^{\lceil (n-m-2)/2\rceil} \p_\theta^{(n-m-2)}(D_{n-m-2,0})}{\binom{n}{m} \theta^{\lceil (n-m)/2\rceil} \p_\theta^{(n-m)}(D_{n-m,0})}\\
		={}& \frac{\binom{n-2}{m} 2 \binom{n-m-2}{3} (n-m-6)!!}{\binom{n}{m} 2 \binom{n-m}{3} (n-m-4)!!}
		= \frac{n-m-3}{n(n-1)}
	\end{align*}
	and
	\begin{align*}
		\lim_{\theta \to \infty} \frac{\p_\theta^{(n-2)}(D_{n-2,m+2})}{\theta^2 \p_\theta^{(n)}(D_{n,m})}
		= \frac{(n-m-3)(n-m-5)}{n(n-1)(m+2)(m+1)}.
	\end{align*}
	Again, it remains to combine these limits with Proposition~\ref{prop:exact} and to simplify the resulting expressions.
\end{proof}

\begin{proof}[Proof of Proposition~\ref{prop:condTheta0}]
	For the time being, assume~$m \lne n-2$. We have $$\lim_{\theta \to 0} \frac{1}{\theta} \sum_{k=0}^{l} \frac{(-\theta)^k \theta^{(l-k)}}{k!(l-k)!} = \frac{1}{l} \1_{l \neq 1}, \quad l \in \N,$$ because only the summand corresponding to~$k=0$ matters asymptotically if~$l \neq 1$, while the sum is~$0$ if~$l=1$. Together with Proposition~\ref{prop:pFix} (for~$m=0$, the calculation is trivial), we get
	\begin{align*}
		\lim_{\theta \to 0} \theta \frac{\p_\theta^{(n-2)}(D_{n-2,m-1})}{\p_\theta^{(n)}(D_{n,m})}
		={}& \lim_{\theta \to 0}\frac{m(\theta+n-1)(\theta+n-2) \sum_{k=0}^{n-m-1} \frac{(-\theta)^k \theta^{(n-m-1-k)}}{k! (n-m-1-k)!}}{n(n-1) \sum_{k=0}^{n-m} \frac{(-\theta)^k \theta^{(n-m-k)}}{k! (n-m-k)!}}\\
		={}& \frac{m(n-2) \frac{1}{n-m-1}}{n \frac{1}{n-m}}.
	\end{align*}
	Similarly, we obtain
	\begin{align*}
		\lim_{\theta \to 0} \frac{\p_\theta^{(n-2)}(D_{n-2,m})}{\p_\theta^{(n)}(D_{n,m})}
		={}& \frac{(n-2) \frac{1}{n-m-2}}{n \frac{1}{n-m}} \1_{m \neq n-3}
	\end{align*}
	as well as
	\begin{align*}		
		\lim_{\theta \to 0} \frac{\p_\theta^{(n-2)}(D_{n-2,m+1})}{\p_\theta^{(n)}(D_{n,m})}
		= 0 \quad \text{and} \quad \lim_{\theta \to 0} \frac{\p_\theta^{(n-2)}(D_{n-2,m+2})}{\p_\theta^{(n)}(D_{n,m})} = 0.
	\end{align*}
	It remains to combine these limits with Proposition~\ref{prop:exact} and to simplify the resulting expressions.

	In the special case~$m=n-2$, we trivially have
	\begin{align*}
		\frac{\p_\theta^{(n-2)}(D_{n-2,(n-2)-1})}{\p_\theta^{(n)}(D_{n,n-2})} = \frac{\p_\theta^{(n-2)}(D_{n-2,(n-2)+1})}{\p_\theta^{(n)}(D_{n,n-2})} = \frac{\p_\theta^{(n-2)}(D_{n-2,(n-2)+2})}{\p_\theta^{(n)}(D_{n,n-2})} = 0, \quad \theta \gne 0.
	\end{align*}
	Using Proposition~\ref{prop:pFix} or elementary combinatorics, we further get
	\begin{align*}
		\lim_{\theta \to 0} \theta \frac{\p_\theta^{(n-2)}(D_{n-2,n-2})}{\p_\theta^{(n)}(D_{n,n-2})}
		={}& \lim_{\theta \to 0} \frac{2(\theta+n-1)(\theta+n-2)}{n(n-1)} = \frac{2(n-2)}{n}.
	\end{align*}
	Together with Proposition~\ref{prop:exact}, this yields the claim for~$m=n-2$.
\end{proof}

\end{document}